\documentclass[numbook, envcountsect, envcountsame, envcountreset, runningheads, final, smallextended]{elsarticle2}

\usepackage{amsmath}
\usepackage{amsthm}
\usepackage{enumitem}
\usepackage[numbers]{natbib}

\newtheorem{thm}{Theorem}[section]
\newtheorem{lem}[thm]{Lemma}
\newtheorem{prop}[thm]{Proposition}
\newtheorem{cor}[thm]{Corollary}

\theoremstyle{definition}
\newtheorem{defn}[thm]{Definition}
\newtheorem{con}[thm]{Conjecture}
\newtheorem{exmp}[thm]{Example}

\DeclareMathOperator{\des}{\text{\rm des}}
\DeclareMathOperator{\susp}{\text{\rm susp}}
\DeclareMathOperator{\rank}{\text{\rm rank}}

\begin{document}
	
\begin{frontmatter}
\title{Intervals of Permutations with a Fixed Number of Descents are Shellable}
\author{Jason P Smith}
\address{University of Strathclyde, Glasgow, UK\\ jason.p.smith@strath.ac.uk}

\begin{abstract}
The set of all permutations, ordered by pattern containment, is a poset. We present an order isomorphism from the poset of permutations with a fixed number of descents to a certain poset of words with subword order. We use this bijection to show that intervals of permutations with a fixed number of descents are shellable, and we present a formula for the M\"obius function of these intervals. We present an alternative proof for a result on the M\"obius function of intervals~$[1,\pi]$ such that $\pi$ has exactly one descent. We prove that if $\pi$ has exactly one descent and avoids 456123 and 356124, then the intervals $[1,\pi]$ have no nontrivial disconnected subintervals; we conjecture that these intervals are shellable.
\end{abstract}

\begin{keyword}
Permutation Poset \sep Shellability \sep M\"obius Function
\MSC 05E45 \sep 05A05
\end{keyword}

\end{frontmatter}

\section{Introduction and Preliminaries}\label{sec:intro}
A \emph{permutation} of length $n$ is an ordering of the integers $1,\ldots,n$, without repetitions. Given two permutations $\sigma$ and $\pi$, we define an occurrence of~$\sigma$ as a \emph{pattern} in $\pi$ to be a subsequence of $\pi$ with the same relative order of elements as in $\sigma$. For example, if $\sigma=213$ and $\pi=23514$ then there are two occurrences of $\sigma$ in $\pi$, as the subsequences $214$ and $314$. The set of all permutations forms a poset $\mathcal{P}$, with a partial ordering defined by~$\sigma\le\pi$ if~$\sigma$ occurs as a pattern in~$\pi$. An \emph{interval} $[\sigma,\pi]$ in $\mathcal{P}$ is a subposet consisting of all permutations $z\in \mathcal{P}$ with~$\sigma\le z\le\pi$. A \emph{chain} in a poset $P$ is a totally ordered subset $\{c_1<\cdots<c_t\}$. For example, $21<2341<24513$ is a chain in $[1,24513]$. The direct sum $\sigma\oplus\pi$ of two permutations $\sigma$ and~$\pi$ is obtained by appending~$\pi$ to $\sigma$ after adding the length of $\sigma$ to each letter of $\pi$. For example, $213\oplus312=213645$. A \emph{descent} occurs at $i$ in a permutation $\pi_1\ldots\pi_n$ if~$\pi_i>\pi_{i+1}$. As an example, $23154$ has descents at 2 and 4. 

If $\sigma\le\pi$, then $\des(\sigma)\le\des(\pi)$. Therefore, if the permutations~$\sigma$ and~$\pi$ both have exactly $k$ descents, then any permutation~$\tau\in[\sigma,\pi]$ also has exactly~$k$ descents. We denote the induced subposet of all permutations with exactly~$k$ descents as $\mathcal{P}_k$. The M\"obius function for a poset is defined recursively as follows:~$\mu(a,b)=0$ if $a\not\le b$, $\mu(a,a)=1$ for all~$a$ and for $a<b$:
$$\mu(a,b)=-\sum_{a\le z<b}\mu(a,z).$$ One of the main goals of this paper is to study the M\"obius function of $P_k$. 

The \emph{interior} of the interval $[\sigma,\pi]$, written $(\sigma,\pi)$, is the set $[\sigma,\pi]-\{\sigma,\pi\}$.  The \emph{order complex} of $[\sigma,\pi]$, written $\Delta(\sigma,\pi)$, is the simplicial complex whose faces are the chains contained in the interior $(\sigma,\pi)$. When we attribute a topological property to an interval we mean the corresponding property of its order complex. We refer the reader to \cite{Wac07} for extensive background on the subject of order complexes.

A simplicial complex is \emph{pure} if all its maximal faces, which are called \emph{facets}, have the same dimension. The order complex of an interval of permutations is always pure. A pure simplicial complex $\Delta$ is \emph{shellable} if its facets can be arranged in linear order $F_1,\ldots,F_t$ in such a way that the subcomplex $\left( \cup_{i=1}^{k-1} \langle F_i \rangle \right) \cap \langle F_k \rangle$ is pure and ($\dim\Delta-1$)-dimensional for ~$2\le k\le t$, where $\langle F \rangle = \{ G : G \subseteq F \}$, that is, $\langle F \rangle$ is the subcomplex generated by $F$. Again we refer the reader to~\cite{Wac07} for extensive background on the subject of shellability.

Let $\mathcal{A}$ be the poset of words on the alphabet of positive integers, with the partial order called \emph{subword order} where $v\le w$, with $w=w_1\dots w_n$, if there is a subsequence $w_{i_1}\ldots w_{i_m}$ in $w$ such that $v=w_{i_1}\ldots w_{i_m}$. For example, $2132\le 212312$ but $2132\not\le 21233$. In \cite{Bjo90} a formula was given for computing the M\"obius function on intervals of $\mathcal{A}$ in polynomial time, and it is shown that all intervals in $\mathcal{A}$ are shellable. In this paper we present an order isomorphism, that is, an order-preserving bijection, from each interval in the permutation posets $\mathcal{P}_k$ to a corresponding interval in $\mathcal{A}$. This allows us to easily compute the M\"obius function of intervals from the posets $\mathcal{P}_k$ and to show that they are~shellable.

The \emph{reduced Betti number} $\tilde{\beta}_k(X)$ of a simplicial complex $X$ is the rank of the~$k$-th reduced homology group of $X$ (for background on the homology of simplicial complexes we refer the reader to \cite{Koz08}). The Philip Hall Theorem and the Euler-Poincar\'e formula, which appear as Proposition 1.2.6 and Theorem~1.2.8 in~\cite{Wac07}, combined state:
\begin{align}\label{for:Euler}\mu(\sigma,\pi)=\tilde{\chi}(\Delta(\sigma,\pi))=\sum_{i=-1}^{\dim\Delta(\sigma,\pi)}(-1)^i \tilde{\beta}_i(\Delta(\sigma,\pi)),\end{align}
where $\tilde{\chi}(\Delta(\sigma,\pi))$ is the reduced Euler characteristic of the order complex\linebreak of $[\sigma,\pi]$.

An important property of simplicial complexes is Cohen-Macaulayness,\linebreak which has its origins in commutative algebra. A simplicial complex $\Delta$ is \emph{Cohen-Macaulay} if $\rank(\tilde{H}_i (\ell k_{\Delta}F))=0$ for all $F\in\Delta$ and $i<\dim\ell k_\Delta F$, where~$\ell k_{\Delta}F$ denotes the \emph{link} of $F$ and $\tilde{H}_i$ denotes the $i$'th reduced homology group. For a full explanation of this definition see \cite[Section 4]{Wac07}. A shellable simplicial complex is Cohen-Macaulay, as observed in \cite{Sta96}. We use this property to compute the homology of intervals from the posets $\mathcal{P}_k$ for any $k\ge0$. 

There is a \textit{generalised subword order}, defined in \cite{SagVat06}, where we take a poset~$P$ and let $P^*$ denote the poset of finite words whose letters are elements of~$P$. If~$u,\, w\in P^*$ then $u\le_{P^*} w$ if there is a subword $w_{i_1}\dots w_{i_{|u|}}$ such that $u_j\le_P w_{i_j}$ for $1\le j\le |u|$. If $P$ is an antichain, then generalised subword order is precisely the subword order. In \cite{SagVat06} a formula was presented for the M\"obius function of words with generalised subword order when $P$ is a chain. That paper also established an order isomorphism between posets of these words and posets of \emph{layered permutations}, that is, permutations that can be expressed as a direct sum of decreasing permutations.  For example, $1\oplus21\oplus321\oplus21=13265487$ is a layered permutation. 

In \cite{McnSag12} a formula was presented for the M\"obius function of words with generalised subword order for any poset $P$, which covers both the words considered in the present paper and in \cite{SagVat06}. In \cite{Bjo80} it was shown that if an interval $\mathcal{I}$ contains a \emph{nontrivial} disconnected subinterval, that is, a disconnected subinterval of rank at least 3, then $\mathcal{I}$ is not shellable. The first major result on the topology of intervals from the poset $\mathcal{P}$ appeared in \cite{McSt13}, where it was shown that if $P$ is a rooted forest, then any interval~$[u,v]$ in~$P^*$ that does not contain a nontrivial disconnected subinterval is shellable. This result was then used to show that intervals of layered permutations that do not contain a nontrivial disconnected subinterval are shellable. Furthermore, it was conjectured that the same applies to the more general class of \emph{separable permutations}, that is, the permutations that avoid $2413$ and $3142$.

In Section \ref{sec:bij} we present a bijection from $\mathcal{P}$ to a subposet of $\mathcal{A}$. We show that when we restrict this bijection to $\mathcal{P}_k$ it is an order isomorphism. This allows us to draw on many useful results that have been proven for subword order, such as the shellability of intervals, and apply these results to permutations. In Section~\ref{sec:mob} we use this order isomorphism to present a formula for the M\"obius function of intervals from the posets $\mathcal{P}_k$. We use this formula to prove a conjecture made in~\cite{Smith13} and to present an alternative, simpler proof of \cite[Theorem~5]{Smith13} on the M\"obius function of intervals $[1,\pi]$ such that $\pi$ has one descent. In Section~\ref{sec:shell} we show that if $\pi$ has exactly one descent and avoids 456123 and 356124, then~$[1,\pi]$ has no nontrivial disconnected subintervals and we conjecture that these intervals are shellable.

\section{Bijection From Permutations to Words}\label{sec:bij}
In this section we present an order isomorphism from the poset $P_k$ of permutations with exactly $k$ descents to a subposet of $\mathcal{A}$. Let $\max(w)$ be the value of the largest letter in the word $w$. We now define the poset of words we consider:

\begin{defn}\label{addCon}
Let $\widehat{\mathcal{A}}$ denote the poset of words with subword order on the alphabet of all positive integers, with the additional conditions that for any~$w\in~\widehat{\mathcal{A}}$:
\begin{description}[font=\normalfont]
\item[AC1:] There is at least one occurrence of each letter $i\in\{1,\ldots,\max(w)\}$.\label{con:0}
\item[AC2:] The rightmost occurrence of each letter $i\in\{1,\ldots,\max(w)-1\}$ is preceded by an occurrence of $i+1$.\label{con:2}
\end{description}
Let $\widehat{\mathcal{A}_k}$ denote the subposet of $\widehat{\mathcal{A}}$ of words $w$ where $\max(w)=k$.
\end{defn}

\begin{exmp}
For example, $231423\in\widehat{\mathcal{A}}$ but $1121343\not\in\widehat{\mathcal{A}}$ because the rightmost occurrence of $2$ does not have a $3$ to its left.
\end{exmp}

The additional conditions in Definition \ref{addCon} are very similar to the definition of a \emph{restricted growth function}, which can be used to encode set partitions, see~\cite{Mil77}. To see the similarity we use the definition of a restricted growth function that appears in Question 106 in \cite[Chapter~1]{Sta96}. A \emph{restricted growth function} is a sequence of the positive integers~$1,\ldots,k$ with each letter occurring at least once and the first occurrence of $i$ appearing before the first occurrence of $i+1$, for~$1\le i \le k-1$. If we consider AC2 reworded as beginning at the right end of the word, and travelling left, then the first occurrence of~$i$ must appear before the last occurrence of~$i+1$. The key difference is that AC2 requires at least one occurrence of $i+1$ after the first~$i$ whereas a restricted growth function requires that all occurrences of $i+1$ are after the first~$i$. As such, it is easy to see that~$\widehat{\mathcal{A}}$ is a larger class than the class of restricted growth functions.

We know that the number of permutations of length $n$ in $\mathcal{P}_k$ is the Eulerian number $A(n,k)$, see \cite{Sta97}. We show that there is a length-preserving bijection from $\widehat{\mathcal{A}_k}$ to~$\mathcal{P}_{k-1}$, which implies the number of words of length $n$ in $\widehat{\mathcal{A}_k}$ is given by the Eulerian number $A(n,k-1)$.

When referring to both words and permutations we often use the notation~$\alpha_i$ to refer to the letter at location $i$ in $\alpha$, and $|\alpha|$ to denote the length of $\alpha$. Given a letter $c$ of the permutation $\pi$, let $d_{\pi}(c)$ be the index of the \emph{run} containing $c$, where a run is a maximal consecutive sequence of increasing letters. Therefore,~$d_{\pi}(c)$ equals the number of descents preceding $c$ in $\pi$, plus 1. For example,~$d_{35241}(5)=1$ and~$d_{35241}(1)=3$. Given a letter $j$ of the word $w$, let $p_w(j)$ be the set of positions of the letter $j$ in $w$, in increasing order. For example, $p_{21232}(2)=\{1,3,5\}$. Now define the following functions:
$$f:\mathcal{P}\rightarrow\widehat{\mathcal{A}}\text{ by }\pi\mapsto d_{\pi}(1)d_{\pi}(2)\ldots d_{\pi}(|\pi|),$$
$$g:\widehat{\mathcal{A}}\rightarrow \mathcal{P}\text{ by }w\mapsto p_w(1),\ldots,p_w(\max(w)).$$

Now consider what these functions are doing. When applying $f$ to $\pi$ we first find the location of 1 in $\pi$ and count the number of preceding descents, which gives the first letter, $d_{\pi}(1)$, of $f(\pi)$. We then repeat this for the letter~2 in~$\pi$ and continue up to $n$. To apply $g$ to $w$ we find the positions of each~1 in~$w$ and   $g(w)$ begins with these positions in increasing order. Then we find the positions of each 2 in $w$ and we continue $g(w)$ with these positions in increasing order. We continue this up to $\max(w)$. For example, if~$\pi= 26 34 15$ then $f(\pi)=312231$, and if $w=214321$ then $g(w)=261543$.

 Before proceeding we define a term used for both permutations and words:
\begin{defn}
Consider two elements $a\le b$ of a poset of either words or permutations. An \emph{embedding} of $a$ in $b$ is a sequence $\eta$ of length~$|b|$ such that the nonzero positions in $\eta$ are the positions of an occurrence of $a$ in $b$ and removal of all the zeros from $\eta$ results in~$a$.
\end{defn}
\begin{exmp}\label{ex:em}
The embeddings of the word 2121 in 211221 are 210201, 210021, 201201 and 201021. The embeddings of the permutation 213 in $142356$ are $021030,\, 020130,\, 021003$ and $020103$.
\end{exmp}

We now show that $f$ and $g$ are inverses of each other. First we show that $f$ and $g$ link the number of descents of permutations in $\mathcal{P}$ and the largest letter of words in $\widehat{\mathcal{A}}$.

\begin{lem}\label{lem:im}
Let $f$ and $g$ be defined as above.
\begin{enumerate}
\item If $\pi\in \mathcal{P}_k$, then $f(\pi)\in\widehat{\mathcal{A}_{k+1}}$.\label{lem:im:1}
\item If $w\in\widehat{\mathcal{A}_{k+1}}$, then $g(w)\in \mathcal{P}_k$.\label{lem:im:2}
\end{enumerate}
\begin{proof}
For \eqref{lem:im:1}, consider $\pi\in \mathcal{P}_k$ and let $w=f(\pi)$. It is clear that~$w$ is a word and that $d_\pi(\pi_n)=k+1$. Also there must be an occurrence of all the letters~$1,\ldots,k$ because for each $i$ the letter at the location of the $i$-th descent maps to $i$. All that remains to be shown is that $w$ satisfies AC2 in Definition~\ref{addCon}. Let $w_t=i$ be the rightmost occurrence of the letter $i$. This implies the letter $t$ at position $j$ in $\pi$ is the rightmost letter that is preceded by exactly $i$ descents, and hence a descent occurs directly after $\pi_j$. Thus the letter $\pi_{j+1}$ is mapped to $i+1$. Since $\pi_{j+1}<\pi_j$, the letter $\pi_{j+1}$ is mapped to an earlier location in~$w$ than $\pi_j$. Therefore, $w_t=i$ is preceded by an occurrence of $i+1$. Since the argument holds for all $i$, this proves \eqref{lem:im:1}.

For \eqref{lem:im:2}, we need to show there are $k$ descents in $g(w)$. By AC2  in Definition~\ref{addCon}, the largest letter in $p_w(t)$ must have a greater value than the smallest letter in $p_w(t+1)$. Therefore, for each $t$ there is a descent between~$p_w(t)$ and~$p_w(t+1)$ in $g(w)$. Since each $p_w(j)$ is increasing, these $k$ are the only descents.

\end{proof}
\end{lem}

\begin{lem}\label{lem:bij}
The map $f$ is a bijection with inverse $g$.
\begin{proof}
We prove this by showing that $f g=id_{\widehat{\mathcal{A}}}$ and $g f=id_{\mathcal{P}}$.

First consider $w\in\widehat{\mathcal{A}}$ and $v=f(g(w))$. If $w_i=t$ then $d_{g(w)}(i)=t$,\linebreak since~$i\in p_w(t)$ and thus in the $t$-th run of $g(w)$. Since we know\linebreak that $w_i=t=d_{g(w)}(i)=v_i$ for all $i$, we conclude $w=v$.

Now consider $g(f(\pi))$ such that $\pi\in\mathcal{P}_k$, and let $\pi_t\ldots\pi_{t+\lambda}$ be the $j$-th run of $\pi$ for some $j\in\{1,\ldots,k+1\}$. Each $\pi_\ell$, where $\ell\in\{t,\ldots,t+\lambda\}$, is mapped to the letter $j$ in~$f(\pi)$, and these are the only letters mapped to $j$. In turn only those letters are mapped into $p_{f(\pi)}(j)$. Since each segment is listed in increasing order, and this holds for all $j$, we have $g(f(\pi))=\pi$.
\end{proof}
\end{lem}

So $f$ is a bijection from $\mathcal{P}$ to $\widehat{\mathcal{A}}$. Finally we need to see if this bijection is order-preserving. This is not true in general. For example, consider the permutations $132\le2143$: Applying $f$ yields $f(132)=121\not\le2132=f(2143)$.

Consider the functions $f_k$ obtained by restricting $f$ to $\mathcal{P}_k$ and $g_k$ obtained by restricting $g$ to $\widehat{\mathcal{A}_{k}}$. We know by Lemma \ref{lem:im} that the image of $f_k$ is $\widehat{\mathcal{A}_{k+1}}$ and the image of~$g_{k+1}$ is $\mathcal{P}_k$. Combining this with Lemma \ref{lem:bij} implies $f_k$ is a bijection. We now show that $f_k$ and $g_k$ are order-preserving:

\begin{thm}\label{thm:op}
The bijection $f_k$ is an order isomorphism.
\begin{proof}
Consider two permutations $\sigma,\pi\in \mathcal{P}_k$ with $\sigma\le\pi$. Since $\sigma$ and $\pi$ have the same number of descents, thus the same number of runs, for any occurrence of~$\sigma$ in $\pi$ the $t$-th run of $\sigma$ must occur in the $t$-th run of~$\pi$. If $\pi_{k_1}\ldots\pi_{k_m}$  is an occurrence of $\sigma$ in $\pi$, then $d_\pi(\pi_{k_1})\ldots d_\pi(\pi_{k_m})=d_\sigma(\sigma_1)\ldots d_\sigma(\sigma_m)$. Let $\pi_{t_1}\ldots\pi_{t_m}$ be the reordering of $\pi_{k_1}\ldots\pi_{k_m}$ in increasing order, then $d_\pi(\pi_{t_1})\ldots d_\pi(\pi_{t_m})$ occurs in $f_k(\pi)$ and is equal to $f_k(\sigma)$. Therefore, $f_k(\sigma)\le f_k(\pi)$.

Now consider two words $v,w\in\widehat{\mathcal{A}_k}$ with $v\le w$. Let $\eta$ be an embedding of~$v$ in~$w$, and let $\widehat{g_k}(\eta)=p_{\eta}(1)\ldots p_{\eta}(k+1)$. It is easy to see that~$p_{\eta}(t) \subseteq~p_{w}(t)$, which implies $\widehat{g_k}(\eta) \le g_k(w)$. Also $\widehat{g_k}(\eta)$ is an occurrence of~$g_k(v)$,\linebreak so~$g_k(v)\le~g_k(w)$.
\end{proof}
\end{thm}

Hence we have an order isomorphism between $\mathcal{P}_k$ and $\widehat{\mathcal{A}_{k+1}}$. One of our key results is the following corollary, which follows directly from \cite[Theorem 3]{Bjo90} and Theorem \ref{thm:op}:

\begin{cor}\label{cor:shell}
Any interval $[\sigma,\pi]$, where $\sigma$ and $\pi$ are permutations with the same number of descents, is dual CL-shellable.
\end{cor}

Note that CL-shellability implies shellability, so a poset that is dual CL-shellable is shellable. For a good survey of the implications of different types of shellability we refer the reader to \cite[Section~4.1]{Wac07}.

We can also consider $f$ as a map to the poset of words on generalised subword order, where the underlying poset is the chain of positive integers. In this case~$f$ is order-preserving, but $g$ is not. For example, $211\le212$ but $g(211)=231\not\le213=g(212)$.

It is known that a shellable complex has the homotopy type of a wedge of spheres. Therefore, Corollary \ref{cor:shell} gives the following result:

\begin{cor}
If $\sigma$ and $\pi$ are permutations with the same number of descents, then $\Delta(\sigma,\pi)$ is homotopy equivalent to a wedge of $|\mu(\sigma,\pi)|$ spheres of dimension~$\dim\Delta(\sigma,\pi)=|\pi|-|\sigma|-1$.
\end{cor}

\section{Computing the M\"obius function}\label{sec:mob}
We can use Theorem \ref{thm:op} along with \cite[Theorem 1]{Bjo90}, which also appears as~\cite[Theorem 2.1]{SagVat06}, to compute the M\"obius function of any interval in $\mathcal{P}$ between permutations with the same number of descents. To do this we first need to define what a normal embedding is in the case of permutations. The definition we use is induced by the definition of a normal embedding in \cite{SagVat06} after applying the bijection from Theorem \ref{thm:op}:

\begin{defn}\label{def:emb}
An \emph{adjacency} in a permutation is a sequence of consecutively valued letters in increasing consecutive order. The \textit{tail} of an adjacency is all but the first letter of the adjacency. An embedding $\eta$ of $\sigma$ in $\pi$ is \textit{normal} if~$\eta_i$ is nonzero for each letter $\pi_i$ in the tail of an adjacency. We use the notation from~\cite{Bjo90} and denote the number of normal embeddings of $\sigma$ in $\pi$ as $\binom{\pi}{\sigma}_n$.
\end{defn}
There is an analogous decreasing adjacency, but we are only interested in increasing adjacencies.
\begin{exmp}
As in Example \ref{ex:em} consider 213 and 142356.  The adjacencies in~$142356$ are 23 and 56 so the tails of the adjacencies are $3$ and $6$. Hence the only normal embedding is $020103$ and therefore $\binom{142356}{213}_n=1$.
\end{exmp}

We use this definition to state the following result:

\begin{prop}\label{prop:mob}
If $\sigma$ and $\pi$ are permutations with the same number of descents, then
$$\mu(\sigma,\pi)=(-1)^{|\pi|-|\sigma|}\binom{\pi}{\sigma}_n.$$
\begin{proof}
This follows directly from Theorem \ref{thm:op} and  \cite[Theorem 1]{Bjo90}.
\end{proof}
\end{prop}

In \cite{Bjo90} it was shown that $\binom{\pi}{\sigma}_n$ can be computed in polynomial time.

In Section~\ref{sub:mob} we use Proposition \ref{prop:mob} to give a simpler proof of a result which appears in \cite{Smith13} and prove a conjecture from the same paper. First we present two corollaries:

\begin{cor}\label{cor:a}
Consider $\sigma,\,\pi\in \mathcal{P}_k$. Let $t$ be the total number of letters in all the tails of all the adjacencies in $\pi$. If $t>|\sigma|$, then $\mu(\sigma,\pi)=0$.
\end{cor}

This result doesn't hold if we remove the restriction on the number of descents. For example, consider $\sigma=213$ and $\pi=569341278$, which have one and two descents, respectively. The total number of letters in all the tails of~$569341278$ is $t=4$ and $|\sigma|=3$, but $\mu(312,6745123)=1\not=0$.

Corollary \ref{cor:a} is another part of the answer to a question posed in \cite{BJJS11} asking when is $\mu(\sigma,\pi)=0$. Whilst we cannot yet give a simple definitive answer to this question, there are results which present several classes of intervals with a zero M\"obius function, such as results in \cite{BJJS11}, \cite{Smith13} and~\cite{SteTen10}.

 A result in \cite{BJJS11} showed that if $\sigma$ and $\pi$ are separable permutations,\linebreak then $|\mu(\sigma,\pi)|$ is at most the number of occurrences of $\sigma$ in $\pi$. Proposition \ref{prop:mob} implies this is also the case if we fix the number of descents, since an embedding corresponds to a unique occurrence.

\begin{cor}\label{cor:bdd}
If $\sigma$ and $\pi$ have the same number of descents, then $|\mu(\sigma,\pi)|$ is at most the number of occurrences of $\sigma$ in $\pi$.
\end{cor}

\subsection{M\"obius Function of Permutations With at Most One Descent}\label{sub:mob}
Proposition \ref{prop:mob} allows us to compute the M\"obius function of an interval between two permutations with the same number of descents, but says nothing about intervals between permutations with different number of descents. Now we consider the intervals $[\textbf{1},\pi]$, where $\pi\in\mathcal{P}_1$ and $\textbf{1}$ denotes the permutation~1. In particular we present an alternative proof, which is both shorter and simpler than the original, of \cite[Theorem 5]{Smith13}. We begin with a useful lemma which gives a formula for $\mu(\textbf{1},\pi)$ for every permutation $\pi$ with one descent.

\begin{lem}\label{1eq2}
If $\pi$ has exactly one descent, then $\mu(\textbf{1},\pi)=-\mu(21,\pi)$.
\end{lem}
Lemma \ref{1eq2} can be proved directly by considering the effect the removal of the increasing permutations has on the M\"obius function. However, it also follows from Theorem \ref{prop:hom}, so we omit the proof here.

We now present the alternative proof of \cite[Theorem 5]{Smith13}. As in \cite{Smith13}, we use the notation $\mu(\pi):=\mu(\textbf{1},\pi)$. A \emph{triple adjacency} indicates an adjacency of three letters, for example 234 in 52341. We use the notation \emph{adjacency pair} to denote an adjacency of length 2. The value and position of an adjacency pair are given by the value and position of the first letter of the adjacency pair.  We denote the two permutations of length $n$ that have one descent and no adjacencies as~$M_n=246\ldots135\ldots$ and $W_n=135\ldots246\ldots$. For example, $M_6=246135$ and~$W_5=13524$.

As observed in \cite{Smith13}, in Theorem~\ref{thm:main} any overlap of cases agree in value. For example, if $\pi$ contains the triple adjacency 234, then equivalently $\pi$ contains the two adjacency pairs 23 and 34, the first of which has lower value;  both cases imply $\mu(\pi)=0$.
\begin{thm}\label{thm:main}
Given a permutation $\pi$ of length $n>2$, with exactly one descent, the value of $\mu(\pi)$ can be computed from the number and positions of adjacencies in $\pi$, as follows:
\begin{enumerate}
\item If $\mu(\pi)\not=0$, then $\mu(\pi)$ is positive if and only if $n$ is odd.\label{main7}
\item If $\pi$ begins with $12$ or ends in  $(n-1)n$, then $\mu(\pi) = 0$.\label{main1}
\item If $\pi$ has a triple adjacency, then $\mu(\pi) = 0$.\label{main2}
\item If $\pi$ has more than two adjacency pairs, then $\mu(\pi) = 0$.\label{main3}
\item If $\pi$ has exactly two adjacency pairs, then:\label{main4}
\begin{enumerate}
\item If the first adjacency pair has greater value than the second,\linebreak then $|\mu(\pi)|=~1$,\label{main4a}
\item If the first adjacency pair has lower value than the second,\linebreak then~$\mu(\pi)=~0$.\label{main4b}
\end{enumerate}
\item If $\pi$ has exactly one adjacency pair, at position $i\in\{1,\ldots,n-1\}$, and the descent is at position $d$, then:\label{main5}
\begin{enumerate}
\item If $i<d$ and $\pi_1\not=1$, then $|\mu(\pi)|=i$,\label{main5a}
\item If $i<d$ and $\pi_1=1$, then $|\mu(\pi)|=i-1$,\label{main5b}
\item If $i>d$ and $\pi_{n}\not=n$, then $|\mu(\pi)|=n-i$,\label{main5c}
\item If $i>d$ and $\pi_{n}=n$, then $|\mu(\pi)|=n-i-1$.\label{main5d}
\end{enumerate}
\item If $\pi$ has no adjacencies, then:\label{main6}
\begin{enumerate}
\item If $n$ is even and $\pi_1=1$, so $\pi=W_{n}$, then $\mu(\pi)=-\dbinom{\frac{n}{2}}{2}$,\label{main6a}
\item If $n$ is even and $\pi_1=2$, so $\pi=M_{n}$, then $\mu(\pi)=-\dbinom{\frac{n}{2}+1}{2}$,\label{main6b}
\item If $n$ is odd, then $\mu(\pi)=\dbinom{\frac{n+1}{2}}{2}$.\label{main6c}
\end{enumerate}
\end{enumerate}
\begin{proof}
By Lemma \ref{1eq2}, we know that $\mu(\pi)=-\mu(21,\pi)$. We can use Proposition~\ref{prop:mob} to compute $\mu(21,\pi)$, which implies the sign of $\mu(21,\pi)$ is given\linebreak by $(-1)^{|\pi|-2}$. Therefore, $\mu(21,\pi)$ is positive if and only if $n$ is even, combining~this with $\mu(\pi)=-\mu(21,\pi)$ gives part \ref{main7}. 

We need to show that the absolute value of $\mu(21,\pi)$, which equals the number of normal embeddings, agrees with each of the cases in the theorem. We refer to the permutation 21 as $\sigma$, to avoid confusion between letters and permutations.

Case \ref{main1}: If $\pi$ begins with 12, then we must embed the 2 of $\sigma$ as the 2 in $\pi$. However, there is no letter after the descent of value less than 2, so we cannot embed the~1 of~$\sigma$ anywhere. Similarly, if $\pi$ ends in $(n-1)n$, then we must embed the 1 of~$\sigma$ as $n$ in~$\pi$. However, this leaves no valid position to embed the 2 of $\sigma$. Therefore, there are no normal embeddings of $\sigma$ in $\pi$. 

Case \ref{main2}: If $\pi$ has a triple adjacency at $\pi_i\pi_{i+1}\pi_{i+2}$, then any normal embedding of $\sigma$ in $\pi$ must be non-zero for $\pi_{i+1}\pi_{i+2}$. Therefore, $\sigma$ must contain 12, which~21 does not. So there are no normal embeddings of $\sigma$ in $\pi$. 

Case \ref{main3} follows directly from Corollary \ref{cor:a}.

Case \ref{main4}: When there are two adjacency pairs, at locations $k$ and $j$, there is only one embedding that might be normal, namely $\eta=\ldots0\eta_{k+1}0\ldots0\eta_{j+1}0\ldots$. If~$\pi_k>\pi_j$, then we can set $\eta_{k+1}=2$ and $\eta_{j+1}=1$. Therefore, there is one normal embedding of $\sigma$ in $\pi$. If $\pi_k<\pi_j$, then there is no way to make $\eta$ an embedding of $\sigma$. Therefore, there are no normal embeddings of $\sigma$ in $\pi$.

Case \ref{main5}: In these cases we must embed one of the letters of 21 in the adjacency pair and can choose an appropriate place for the other letter. Denote the locations of the descent and adjacency pair as $d$ and $i$, respectively. If $i<d$, then an embedding~$\eta$ of $\sigma$ in $\pi$ must have $\eta_{i+1}=2$ and we can then embed the~1 from~$\sigma$ in any of the letters after the descent that have value less than $\pi_i$. Since the rest of $\pi$ follows the same alternating pattern, because there are no more adjacencies, it is easy to see that this gives the desired results. The argument is analogous if $i>d$.

Case \ref{main6}: Since there are no adjacencies in $\pi$, any embedding is normal. Therefore, we need only count the number of embeddings. First consider the case when~$n$ is even and $\pi_1=1$. If we embed the letter 2 of $\sigma$ in locations $1,\, 2,\, \ldots \frac{n}{2}$ and then count where we can embed the letter 1, then we get the following sequence $0,1,2,\ldots,\frac{n}{2}-1$. Summing the sequence implies $\binom{\pi}{21}_n=\binom{\frac{n}{2}}{2}$. Repeating this for each case gives the desired results.
\end{proof}
\end{thm}

We can also use Proposition \ref{prop:mob} to prove one of the conjectures presented in~\cite{Smith13}. In Proposition \ref{prop:bottom} we count the number of adjacency pairs, so a triple adjacency counts as two adjacency pairs and a length $k$ adjacency counts as~$k-1$ adjacency pairs. We say that two permutations with exactly one descent are \emph{related} if they have the letter 1 on the same side of the descent. Let $\lfloor x\rfloor$ denote the \emph{floor} of $x$, that is, the largest integer not greater than $x$. 

\begin{lem}\label{lem:big}
Let $\sigma$ be a permutation of length $m$ with exactly one descent and~$i$ adjacency pairs. In $\sigma$ the letter $m$ occurs on the same side of the descent as the letter~1 if and only if $m-i$ is odd.
\begin{proof}
If $\sigma$ begins with the letter 1, then let $\tau=W_m$, otherwise let $\tau=M_m$. We can build $\sigma$ from $\tau$ by going through each letter $k\in\{2,\ldots,m\}$ in $\tau$. If $k$ is not on the same side of the descent in $\tau$ as $k$ is in $\sigma$, then move $k$ to the opposing side of the descent, in the unique way that does not create a new descent. 

We consider three cases that occur when moving a letter $k\in\{2,\ldots,m-~1\}$. If~$k$  is not part of an adjacency pair, then moving it creates two new adjacency pairs $(k-1)k$ and~$k(k+1)$. If $k$ is part of one adjacency pair, then moving it destroys one adjacency pair but creates another. If $k$ is part of two adjacency pairs, then moving it destroys both adjacency pairs. If~$k=m$, then moving it either creates or destroys the adjacency pair~$(m-1)m$. Therefore, each move of a letter~$k$ changes the number of adjacency pairs by $-2$, $0$ or $2$ for all~$k\in~\{2,\ldots,m-1\}$ and by $1$ or $-1$ if $k=m$.

If $m$ is odd, then $1$ and $m$ are on the same side of the descent in $\tau$. If~$m$ is not moved whilst building $\sigma$ from $\tau$, then $m-i$ must be odd and $m$ must be on the same side of the descent as $1$ in $\sigma$. If $m$ is moved, then it is on the opposite side of the descent and $m-i$ is even. The argument is analogous if $m$ is even.
\end{proof}
\end{lem}

\begin{prop}\label{prop:bottom}
Given a permutation $\sigma\in \mathcal{P}_1$ of length $m$, let $i$ be the number of adjacency pairs in $\sigma$.  If $\sigma\le\pi$ such that $\pi\in\{M_n,W_n\}$, then:
\begin{equation*}\mu(\sigma,\pi)=(-1)^{n-m}\dbinom{\lfloor\frac{n+m-i-a}{2}\rfloor}{m},\end{equation*}
where
$a=\begin{cases}
0,&\mbox{ if } \sigma \mbox{ and } \pi \mbox{ are related}\\
1,&\mbox{ otherwise }
\end{cases}.$
\begin{proof}
Since both $\sigma$ and $\pi$ have exactly one descent, we can apply Proposition~\ref{prop:mob}. The sign part of the result follows immediately. Since $\pi$ has no adjacencies, any embedding of $\sigma$ in $\pi$ is normal, hence we need only count the number of embeddings. To do this we find it simpler to consider $f(\sigma)$ and $f(\pi)$, which are binary strings. Note that we consider an occurrence of a substring to occur in consecutive positions. For example, $101$ has an occurrence of $10$, but no occurrence of $11$. We consider the different cases depending on whether $n$ is odd or even, and whether $\sigma$ and $\pi$ are related.

First consider the case when $\sigma$ and $\pi$ are related and $n$ is even. Suppose~$\pi=~M_n$, then $f(\pi)=1010\ldots$ and can be split into $n/2$ blocks, each consisting of a single 10. We can choose to embed a 10 from $f(\sigma)$ in either a single block of~$f(\pi)$ or two separate blocks. For any other letter of $f(\sigma)$ we choose a single block of $f(\pi)$ in which to embed it. Thus, once we decide which~10s of~$f(\sigma)$ to embed in single blocks of $f(\pi)$, all we need to do to determine an embedding is to pick a subset of blocks of $f(\pi)$.

Suppose we embed none of the 10s of $f(\sigma)$ in a single block of $f(\pi)$. We need to pick $m$ of the $n/2$ blocks of $f(\pi)$, to embed one letter of $f(\sigma)$ in each of the selected blocks, which can be done in $\binom{\frac{n}{2}}{m}$ ways. Suppose we select $r$ of the~10s in~$f(\sigma)$ to embed in a single block. We need to choose $m-r$ of the~10s in~$f(\pi)$  in which to embed the parts of $f(\sigma)$. Thus, we need to pick a total of~$m$ objects, some of them blocks of $f(\pi)$ to embed in and some of them~10s in~$f(\sigma)$ to embed in a single block of $f(\pi)$. An occurrence of 10 in~$f(\sigma)$ corresponds to a letter in~$\sigma$ that is after the descent and not the start of an adjacency pair, and there are~$\lfloor\frac{m-i}{2}\rfloor$ such letters. Therefore, we have~$\binom{\lfloor\frac{n}{2}\rfloor+\lfloor\frac{m-i}{2}\rfloor}{m}$ embeddings, and because~$n$ is even this gives the desired result. If $\pi=W_n$, $n$ is even and~$\sigma$ and $\pi$ are related, then the proof is analogous to when $\pi=M_n$, but considering substrings~01 instead of 10.

Now consider the case when $n$ is odd and $\sigma$ and $\pi$ are related. By Lemma~\ref{lem:big} we know that the largest letters in $\sigma$ and~$\pi$ are on same sides of the descent if and only if $m-i$ is odd. Therefore, if $m-i$ is even, then we cannot embed anything in the final letter of $\pi$; thus this case is equivalent to when $n$ is even and $\sigma$ and~$\pi$ are related. If~$m-i$ is odd, then we can embed a letter of $\sigma$ in the largest letter of~$\pi$; thus we have~$\frac{n+1}{2}$ blocks of $f(\pi)$ to embed in. The remaining argument is analogous to when $n$ is even, using the fact that as $n$ and $m-i$ are odd $\lfloor\frac{n+1}{2}\rfloor+\lfloor\frac{m-i}{2}\rfloor=\lfloor\frac{n+m-i}{2}\rfloor$.

Finally consider the cases when $\sigma$ and $\pi$ are not related. In these cases we cannot embed anything in the first letter of $\pi$. Therefore, we can remove the first letter from $\pi$ without changing the number of embeddings. So these cases are equivalent to when $\sigma$ and $\pi$ are related and $\pi$ is of length $n-1$, the latter point accounting for the $-a$ in the equation.
\end{proof}
\end{prop}

\section{Intervals of $[\textbf{1},\pi]$ Where $\pi$ Has One Descent}\label{sec:shell}
We have shown that intervals between two permutations with the same number of descents are shellable. Now we consider intervals of the form $[\textbf{1},\pi]$ such that~$\pi\in~\mathcal{P}_1$. First we present a useful tool called the Quillen Fiber Lemma, which can be found as Theorem~15.28 in \cite{Koz08}. Define the \emph{upper ideal} as~$Q_{\ge x}:=\{y\in Q\,:\, y\ge x\}$.

\begin{prop}\label{thm:Quil}(Quillen Fiber Lemma)
Let $\phi:P\rightarrow Q$ be an order-preserving map between posets such that for any $x\in Q$ the complex $\Delta(\phi^{-1}(Q_{\ge x}))$\linebreak is contractible. Then the induced map between simplicial complexes\linebreak $\Delta(\phi):\Delta(P)\rightarrow\Delta(Q)$ is a homotopy equivalence.
\end{prop}

Note that the order complex of an upper ideal $Q_{\ge x}$ is always contractible to the point $x$. Now we consider the homology of the order complexes of intervals~$[\textbf{1},\pi]$ such that $\pi\in \mathcal{P}_1$.

\begin{thm}\label{prop:hom}
If $\pi\in \mathcal{P}_1$, then the order complex $\Delta(\textbf{1},\pi)$ is homotopy equivalent to a suspension of $\Delta(21,\pi)$. Therefore, the reduced Betti numbers of $\Delta(\textbf{1},\pi)$ are $\tilde{\beta}_n(\Delta(\textbf{1},\pi))=\tilde{\beta}_{n-1}(\Delta(21,\pi))$, for $n>0$, and $\tilde{\beta}_0(\Delta(\textbf{1},\pi))=0$.
\begin{proof}
Let $X=(\textbf{1},\pi)$ and $A=X\setminus[123,\textbf{k}]$, where $\textbf{k}=1\ldots k$ is the largest increasing permutation that occurs in $\pi$. The only permutations in $A$ not in~$(21,\pi)$ are 21 and 12. The permutations 21 and 12 occur as a pattern in every permutation in $(21,\pi)$. Therefore, in the order complex of $A$ each of the points associated to 12 and 21 is the apex of a cone over $\Delta(21,\pi)$, so~$\Delta(A)$ is a suspension of $\Delta(21,\pi)$.

 We use the Quillen Lemma to show that~$\Delta(X)$ is homotopically equivalent to $\Delta(A)$. Consider the map $f:~X\rightarrow A$ defined by:

$$\displaystyle f(\sigma)=\begin{cases}
12, \text{ if } \sigma\in P_0\\
\sigma, \text{ if } \sigma\in P_1
\end{cases}.$$ 
This map is order-preserving and $f^{-1}(A_{\ge a})=X_{\ge a}$ which is an upper ideal, thus~$\Delta(f^{-1}(A_{\ge a}))$ is contractible. Therefore, by the Quillen Fiber Lemma,~$f$ induces a homotopy equivalence between $\Delta(X)$ and~$\Delta(A)$. Thus, $\Delta(X)$ is homotopically equivalent to a suspension of~$\Delta(21,\pi)$. The result on the reduced Betti numbers then follows directly from the property of the suspension that~$\tilde{H}_{n+1}(\susp X)=\tilde{H}_{n}(X)$.~\end{proof}
\end{thm}

It is not true that all intervals $[\textbf{1},\pi]$, $\pi\in P_1$, are shellable, as can be seen by the following example:

\begin{exmp}\label{ex:noshell}
Consider the permutations $456123$ and $356124$. In the interval~$[\textbf{1},456123]$ the subinterval $[123,456123]$ is disconnected and of rank 3, which implies $[\textbf{1},456123]$ is not shellable. Similarly in $[\textbf{1},356124]$ the subinterval~$[123,356124]$ is disconnected and of rank 3. Consequently, if a permutation~$\pi\in \mathcal{P}_1$ contains 456123 or 356124 the interval $[\textbf{1},\pi]$  is not shellable.
\end{exmp}

Whilst it is not true that the intervals $[\textbf{1},\pi]$ are all shellable, we conjecture that containing 456123 or 356124 are the only obstructions to shellability for the intervals $[\textbf{1},\pi]$ when $\pi\in\mathcal{P}_1$.

\begin{con}
If $\pi\in P_1$ and $\pi$ avoids 456123 and 356124, then the interval~$[\textbf{1},\pi]$ is shellable.
\end{con}

We have been unable to prove this conjecture, but we show that these intervals have no nontrivial disconnected subintervals. We prove this below, but first we need a result from \cite{McSt13} and the following definition:
\begin{defn}
Let $\eta$ be an embedding of $\sigma$ in $\pi$. The \emph{zero set} of $\eta$, which we denote $Z_{\eta}$, is the set $\{i:\eta_i=0\}$. The zero set $Z_E$ of a set of embeddings $E$ is the union of the zero sets of all the embeddings in the set $E$.
\end{defn}
\begin{exmp}
Let $\sigma=213$ and $\pi=245136$. Consider the following embeddings of $\sigma$ in $\pi$: $\eta_1=200130$, $\eta_2=200103$ and $\eta_3=020103$. These embeddings have zero sets $Z_{\eta_1}=\{2,3,6\}$, $Z_{\eta_2}=\{2,3,5\}$ and $Z_{\eta_3}=\{1,3,5\}$. Therefore, the set $\{\eta_1,\eta_2,\eta_3\}$ has zero set $\{1,2,3,5,6\}$.
\end{exmp}

\begin{lem}\label{prop53}(see \cite[Proposition 5.3]{McSt13})
Consider two permutations $\sigma<\pi$ such that $|\pi|-|\sigma|\ge 3$. The interval $[\sigma,\pi]$ is not disconnected if the embeddings of~$\sigma$ in~$\pi$ cannot be partitioned into two non-empty sets $E_1$ and $E_2$ such that~$Z_{E_1}\cap~Z_{E_2}=\emptyset$.
\end{lem}

\begin{prop}\label{prop:noDis}
If $\pi\in\mathcal{P}_1$ and $\pi$ avoids 456123 and 356124, then the interval~$[\textbf{1},\pi]$ has no disconnected subintervals of rank 3 or more.
\begin{proof}
By Corollary \ref{cor:shell} we know that intervals between two permutations in~$\mathcal{P}_1$ are shellable, hence have no disconnected subintervals. All that remains is subintervals of the form $[\alpha,\beta]$, of rank 3 or more, with $\alpha\in \mathcal{P}_0$ (so $\alpha$ is an increasing permutation) and $\beta\in \mathcal{P}_1$. We show there is no way to split the embeddings of $\alpha$ in $\beta$ into two sets with disjoint zero sets. To do this we separate the embeddings into three disjoint sets: 
\begin{enumerate}
\item Embeddings with all of $\alpha$ embedded before the descent in $\beta$ constitute the set $E_1$.
\item  Embeddings with all of $\alpha$ embedded after the descent in $\beta$ constitute the set~$E_2$.
\item  Embeddings with part of $\alpha$ embedded before the descent in $\beta$, and part after, constitute the set $E_3$.
\end{enumerate}
Note that each embedding in $E_1$ has zeros in all positions after the descent. Similarly, all embeddings in $E_2$ have zeros in all positions before the descent. Therefore, it is not possible to split $E_1$ or $E_2$ into smaller sets that have disjoint zero sets. Moreover, $E_3$ cannot be split into smaller sets with disjoint zero sets. To see this note that, it is always possible to swap a nonzero letter with a zero letter directly to the right if after the descent, or directly to the left if before the descent. We can use this to build a sequence of embeddings between any two embeddings in $E_3$, where the elements in each adjacent pair in the sequence have only one letter differing in their zero sets. If the zero sets differ by only one element they cannot be disjoint. Since we can build such a sequence between any two embeddings in $E_3$, it is not possible to split $E_3$ into two sets with disjoint zero sets.

Suppose that all three sets are non-empty. Since both $E_1$ and $E_2$ are non-empty, it is not possible to make an embedding that uses all letters from one side of the descent and some letters from the other. This means that each embedding in $E_3$ must have a zero on both sides of the descent. So all embeddings in $E_1$ must be placed in the same set, all embeddings in $E_3$ must be placed in the same set as the embeddings in $E_1$ and all embeddings in $E_2$ must be placed in the same set as the embeddings in $E_3$. So we cannot split the embeddings into two sets with disjoint zero sets.

We now analyse three cases, depending on which of the three sets are empty.

First suppose $E_1$ is empty and that $E_2$ and $E_3$ are non-empty. Consider the embeddings in $E_3$. Unless an embedding embeds all its letters before the descent, and then some after, it has a zero before the descent, so must be put into the same set as $E_2$. Furthermore, as $E_3$ cannot be split into two sets with disjoint zero sets, the only way for $E_2$ and $E_3$ to have disjoint zero sets is if all the embeddings in $E_3$ have no zeros before the descent. We show that the only way such an embedding can exist is if $\beta=\beta_1\beta_2...\beta_d...\beta_i...\beta_n$ with $\beta_i>\beta_d$ and any letter strictly between $\beta_d$ and $\beta_i$ is less than $\beta_d$. Also the number of letters not between~$\beta_d$ and $\beta_i$ must be exactly $|\alpha|$. Therefore, we can embed $\alpha$ as
$$\eta=\alpha_1...\alpha_d0...0\alpha_{d+1}...\alpha_a,$$
such that $\alpha_{d+1}$ is embedded in position $i$. To see this is the only possible embedding suppose there is another embedding $\hat{\eta}\not=\eta$. Since there cannot be a zero before the descent there must be a zero after $\hat{\eta}_i$. This implies it would also be possible to embed the sequence~$\alpha_d...\alpha_a$ after the descent, leaving a zero before the descent, contradicting our requirement for $E_3$. 

If $\eta$ is a valid embedding, then $\beta_{d-2}\beta_{d-1}\beta_d$ must be of one of two forms, either~$c(c+1)(c+2)$ or  $c(c+2)(c+3)$. Otherwise we could build valid embeddings of the form
$$\alpha_1...\alpha_{d-2}00...0\alpha_{d-1}\alpha_d\alpha_{d+1}...\alpha_a,$$
which has a zero before the descent, contradicting our requirement for $E_3$. We also know that there are $|\beta|-|\alpha|\ge3$ letters smaller than~$\beta_d$ that occur after~$\beta_d$. Therefore, the embedding~$\eta$ can only exist if there is an occurrence of either 456123 or 356124 in $\beta$. Since~$\beta$ avoids both these permutations $\eta$ cannot be a valid embedding. So if $E_1$ is empty the embeddings cannot be split into disjoint zero sets.

An analogous argument shows that if $E_2$ is empty, then the embeddings cannot be split into disjoint zero sets.

Now suppose  $E_3$ is empty but $E_1$ and $E_2$ are not. As $E_3$ is empty there can be no increasing sequence of length $|\alpha|$ spread across both sides of the descent. Using this we can repeat the same argument as above showing that~$\beta_{d-2}\beta_{d-1}\beta_d$ must be of one of the forms $c(c+1)(c+2)$ or  $c(c+2)(c+3)$. Therefore, if $\beta$ avoids 456123 and 356124 this case cannot arise.

Therefore, if $\pi\in P_1$ and $\pi$ avoids 456123 and 356124, then for any $\textbf{1}\le\alpha\le\beta\le\pi$ the embeddings of $\alpha$ in $\beta$ cannot be split into two sets with disjoint zero sets. Thus, by Lemma \ref{prop53}, the interval $[\alpha,\beta]$ cannot be disconnected. Therefore,~$[\textbf{1},\pi]$ has no disconnected subintervals of rank 3 or more.
\end{proof}
\end{prop}

\section*{Acknowledgements}
I would like to thank the authors of the GAP system \cite{GAP4} and the authors of the GAP homology package which were very useful in analysing the order complexes studied in Section~\ref{sec:shell}.  I would also like to thank Nik Ru\v{s}kuc for the suggestion to consider permutations with one descent as binary strings, Einar Steingr{\'i}msson for numerous helpful comments, Russ Woodroofe for pointing out an error in an earlier version of the paper and the referees whose comments greatly improved the readability of the paper.

\bibliographystyle{plainnat}
 \newcommand{\noop}[1]{}

\end{document}